\documentclass[a4paper,11pt]{amsart}
\usepackage{amssymb,amsfonts,amsxtra,amscd,mathrsfs}
\usepackage[all]{xy}
\usepackage{fullpage}
\usepackage{graphicx}
\theoremstyle{plain}
\newtheorem{theorem}{Theorem}[section]
\newtheorem{lemma}[theorem]{Lemma}
\newtheorem{cor}[theorem]{Corollary}
\newtheorem{prop}[theorem]{Proposition}
\theoremstyle{definition}
\newtheorem{defi}[theorem]{Definition}

\theoremstyle{remark}
\newtheorem{rem}[theorem]{Remark}
\numberwithin{equation}{section}
\newcommand{\ai}{\ensuremath{A_\infty}}
\newcommand{\forms}[2][\bullet]{\ensuremath{\Omega^{#1}(#2)}}
\newcommand{\drham}[2][\bullet]{\ensuremath{DR^{#1}(#2)}}

\newcommand{\gc}[1][\bullet]{\ensuremath{\mathcal{G}_{#1}}}
\newcommand{\gh}[1][\bullet]{\ensuremath{H_{#1}\mathcal{G}}}
\newcommand{\kgc}[1][\bullet]{\ensuremath{\mathcal{KG}_{#1}}}
\newcommand{\kgh}[1][\bullet]{\ensuremath{H_{#1}\mathcal{KG}}}
\newcommand{\lgc}[1][\bullet]{\ensuremath{\mathcal{LG}_{#1}}}
\newcommand{\lgh}[1][\bullet]{\ensuremath{H_{#1}\mathcal{LG}}}
\newcommand{\ce}[2][\bullet]{\ensuremath{C_{#1}(#2)}}
\newcommand{\hce}[2][\bullet]{\ensuremath{H_{#1}(#2)}}
\newcommand{\rce}[3][\bullet]{\ensuremath{C_{#1}(#2;#3)}}
\newcommand{\hrce}[3][\bullet]{\ensuremath{H_{#1}(#2;#3)}}
\newcommand{\mspc}{\ensuremath{\mathcal{M}_{g,n}}}
\newcommand{\dmcmp}{\ensuremath{\overline{\mathcal{M}}_{g,n}}}
\newcommand{\kcmp}{\ensuremath{\mathcal{K}\overline{\mathcal{M}}_{g,n}}}
\newcommand{\lcmp}{\ensuremath{\mathcal{L}\left[\dmcmp\times\Delta_{n-1}\right]}}
\newcommand{\msup}{\ensuremath{\mathcal{M}^{\Delta}_{\sqcup}}}
\newcommand{\ksup}{\ensuremath{\mathcal{KM}^{\Delta}_{\sqcup}}}
\newcommand{\lsup}{\ensuremath{\mathcal{LM}^{\Delta}_{\sqcup}}}
\newcommand{\gf}{\ensuremath{\mathbb{Q}}}
\newcommand{\dilim}[2]{\ensuremath{\varinjlim_{#1} #2}}
\newcommand{\innprod}{\ensuremath{\langle -,- \rangle}}
\newcommand{\noproof}{\begin{flushright} \ensuremath{\square} \end{flushright}}

\DeclareMathOperator{\Der}{Der}

\begin{document}
\title{Noncommutative geometry and compactifications of the moduli space of curves}
\author{Alastair Hamilton}
\address{University of Connecticut, Mathematics Department, 196 Auditorium Road, Storrs, CT 06269. USA.}
\email{hamilton@math.uconn.edu}
\begin{abstract}
In this paper we show that the homology of a certain natural compactification of the moduli space, introduced by Kontsevich in his study of Witten's conjectures, can be described completely algebraically as the homology of a certain differential graded Lie algebra. This two-parameter family is constructed by using a Lie cobracket on the space of noncommutative 0-forms, a structure which corresponds to pinching simple closed curves on a Riemann surface, to deform the noncommutative symplectic geometry described by Kontsevich in his subsequent papers.
\end{abstract}
\thanks{The work of the author was supported by the Max-Planck-Institut f\"ur Mathematik, Bonn.}
\maketitle
\tableofcontents

\section{Introduction}

\subsection{Background}

Consider the moduli space $\mspc$ of compact Riemann surfaces of genus $g$ with $n$ marked points such that $\chi:=2-2g-n<0$ and $n\geq 1$. In the 1980s it was discovered by the work of mathematicians such as Harer, Mumford, Penner and Thurston that this space admits a description in terms of an orbi-cellular complex where every orbi-cell is indexed by a type of graph called a \emph{ribbon graph}, which lies embedded in the Riemann surface.

In his seminal 92 paper \cite{kontairy}, Kontsevich introduced a certain compactification of the moduli space $\mspc$, which played a crucial role in his proof of Witten's conjectures. An essential point in the proof was that this compactification also admits a description in terms of an orbi-cellular complex.

This compactification was defined as a natural quotient of the Deligne-Mumford compactification $\dmcmp$ by a certain equivalence relation. Although the resulting quotient does not enjoy the good geometric properties of the Deligne-Mumford compactification, in particular it is no longer an orbifold and hence we cannot talk of Poincar\'e duality, it has the advantage that it admits the aforementioned orbi-cellular decomposition and that the tautological classes have a natural description in this framework.

In his subsequent papers \cite{kontsympgeom}, \cite{kontfeynman}; Kontsevich (and later Ginzburg \cite{ginzsympgeom}) developed a framework for studying the symplectic geometry of noncommutative spaces, building on the foundations laid by the work of Connes. In this paper he introduced a certain noncommutative analogue of the Poisson algebra of Hamiltonian vector fields on a symplectic manifold. He demonstrated that the Chevalley-Eilenberg homology of this Lie algebra precisely recovers the homology of the orbi-cellular complex of ribbon graphs and hence the homology of the moduli space $\mspc$.

The goal of this paper is to prove the analogous statement for the compactification of the moduli space introduced by Kontsevich in \cite{kontairy}. This is done by introducing the extra structure of a \emph{Lie bialgebra} on Kontsevich's Lie algebra of noncommutative Hamiltonians. This Lie bialgebra structure is used to construct a two-parameter family of differential graded Lie algebras; one may consider it to be a type of deformation of the original noncommutative symplectic geometry of Kontsevich. It is then shown that the Chevalley-Eilenberg homology of this differential graded Lie algebra recovers the homology of the above compactification.

This leads to the possibility of studying the (co)homological aspects of what is in principle a geometric object, in a purely algebraic manner. In particular, there is a natural way to produce classes in the homology of any differential graded Lie algebra by exponentiating elements in the associated Maurer-Cartan moduli space. In this context the above result provides a way to construct homology classes on this compactification of the moduli space using entirely algebraic data. A description of the resulting algebraic structures, which arise as deformations of $\ai$-structures, will be provided by the author in a subsequent paper.

The algebraic structures which arise in this paper appear in various guises elsewhere in the mathematical literature. Lie bialgebra structures on the space of noncommutative 0-forms appear in the works of Ginzburg and Schedler \cite{ginsch}, \cite{sched}. Algebraic structures related to compactifications of the moduli space are treated by Barannikov in the context of modular operads \cite{baran}. Many relevant ideas appear in the work of Movshev \cite{movshev}. I was also present at a conference at the Institut-Henri-Poincar\'e where I heard Fukaya give a talk \cite{fukaya} on similar algebraic structures and their role in open Gromov-Witten theory. In his talk he explained how he had developed these structures using the framework of string topology provided by Chas-Sullivan \cite{chasul}.

\subsection{Layout of the paper}

The layout of the paper is as follows. In Section \ref{sec_js} we recollect how the theory of Jenkins-Strebel differentials on stable curves leads to certain compactifications of the moduli space that have a natural orbi-cellular structure. In Section \ref{sec_ncgeom} we recall the basic apparatus of noncommutative symplectic geometry introduced in \cite{kontsympgeom} and define a Lie bialgebra structure on the space of noncommutative $0$-forms on a symplectic vector space. In Section \ref{sec_dgla} we use this Lie bialgebra structure to define a two-parameter family of differential graded Lie algebras. In Section \ref{sec_ribgraph} we give a precise formulation of the complex of stable ribbon graphs introduced by Kontsevich in \cite{kontairy} and its relationship with the aforementioned compactifications of the moduli space. In Section \ref{sec_main} we formulate and prove our main theorem which states that the stable homology of the family of differential graded Lie algebras introduced in Section \ref{sec_dgla} recovers exactly the homology of the above compactifications of the moduli space. Throughout the paper we work over the field of rational numbers $\mathbb{Q}$.

\subsection{Jenkins-Strebel theory} \label{sec_js}

Let us recall how the Jenkins-Strebel theory gives us an orbi-cellular decomposition of the moduli space of curves. More detailed expositions can be found in: \cite{looi}, \cite{mondello} and \cite{zvonkine}, on which the following account is based.

Let $R$ be a Riemann surface with $n\geq 1$ marked points and of genus $g>1-\frac{1}{2}n$. A meromorphic section $\beta$ of the tensor square of the holomorphic cotangent bundle of $R$ is called a \emph{quadratic differential}. A \emph{horizontal trajectory} of $\beta$ is a curve on $R$ such that the pullback of the quadratic differential is defined by a \emph{positive} real function.

Given positive real numbers $p_1,\ldots,p_n$; we can invoke the results of Jenkins \cite{jenkins} or Strebel \cite{strebel} which asserts that there exists a \emph{unique} quadratic differential $\beta$ on $R$ such that the following holds:
\begin{enumerate}
\item it has a double pole at each marked point and no other poles,
\item the residue of the pole at the $i$th marked point is $-(\frac{p_i}{2\pi})^2$,
\item the union of all the \emph{closed} horizontal trajectories of $\beta$ is a dense subspace of $R$.
\end{enumerate}
Such a differential is called a \emph{Jenkins-Strebel differential}.

The trajectories of a Jenkins-Strebel differential can be used to decompose the Riemann surface. The union of all the closed horizontal trajectories of our Jenkins-Strebel differential carves out a disconnected open subspace of $R$, whose connected components will be open disks containing one and only one of the marked points. The closed horizontal trajectories surrounding the $i$th marked point will all have length $p_i$ in the metric naturally determined by our quadratic differential.

\begin{figure}[htp]
\centering
\includegraphics{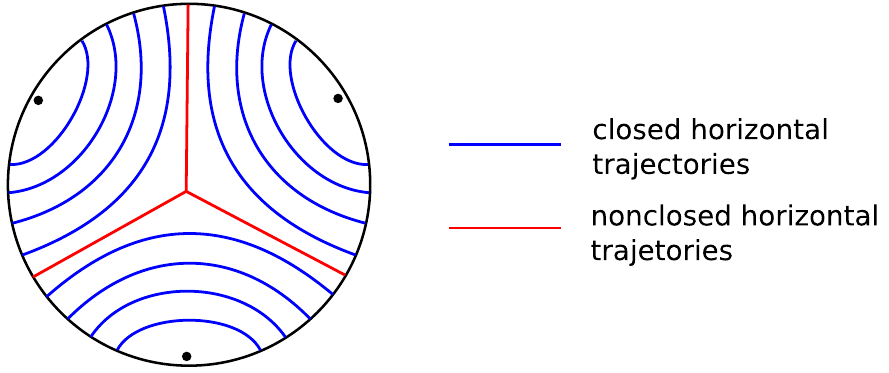}
\caption{Trajectories of a Jenkins-Strebel differential.} \label{fig_trajectories}
\end{figure}

The complement of the subspace carved out by the closed horizontal trajectories is a connected graph which lies embedded in the surface, called the \emph{critical graph} of the Jenkins-Strebel differential. Its edges are the \emph{nonclosed} horizontal trajectories and its vertices are the zeroes of the quadratic differential. A zero of order $n$ gives rise to a vertex of valency $n+2$ and therefore each vertex will be at least trivalent. The embedding of the graph into the surface, together with the natural orientation on the surface, canonically endows every vertex of the graph with a cyclic ordering of the incident half-edges. The resulting graph is called a \emph{ribbon graph}. If we further decorate every edge of the graph by the positive real number corresponding to the length of this edge in the metric determined by the quadratic differential, we obtain the definition of a \emph{metric ribbon graph}.

Hence the Jenkins-Strebel theory provides a way to associate a metric ribbon graph to any Riemann surface. Conversely, given any metric ribbon graph, a standard gluing construction provides a way to reconstruct the corresponding Riemann surface. If we sum the lengths of all the edges of the metric ribbon graph which surround the $i$th marked point, then we recover the positive real number $p_i$, which we call the \emph{perimeter} of the marked point.

This correspondence between Riemann surfaces and metric ribbon graphs leads to the following orbi-cellular decomposition of the decorated moduli space $\mspc\times\Delta^\circ_{n-1}$, due independently to: Harer, Mumford, Penner and Thurston. To any point in $\mspc\times\Delta^\circ_{n-1}$ we can associate a ribbon graph by simply taking the critical graph of the unique Jenkins-Strebel differential on the Riemann surface whose perimeters are prescribed by the coordinate functions on the open simplex $\Delta^\circ_{n-1}$. Now we say that two decorated Riemann surfaces in $\mspc\times\Delta^\circ_{n-1}$ are equivalent if the corresponding ribbon graphs are isomorphic. This equivalence relation partitions the space $\mspc\times\Delta^\circ_{n-1}$ into orbi-cells, each orbi-cell being indexed by a certain ribbon graph.

Since the decorated moduli space $\mspc\times\Delta^\circ_{n-1}$ is not compact, it cannot be an orbi-cellular complex; therefore we compactify it by adding one point. This leads to the following theorem:

\begin{theorem} \label{thm_orbimod}
The one-point compactification of the decorated moduli space $\mspc\times\Delta^\circ_{n-1}$ is an orbi-cellular complex whose orbi-cells are indexed by ribbon graphs (and one 0-cell for the point). An orbi-cell $E$ lies on the boundary of another orbi-cell $E'$ if and only if the ribbon graph corresponding to $E$ is obtained from the ribbon graph corresponding to $E'$ by collapsing edges which are not loops.
\end{theorem}
\noproof

The Jenkins-Strebel theory can also be applied to the Deligne-Mumford compactification of the moduli space. Unfortunately, however, this does not lead to an orbi-cellular decomposition of the Deligne-Mumford compactification. Instead, we obtain an orbi-cellular decomposition of a quotient of the Deligne-Mumford compactification by a certain equivalence relation. The information which is lost corresponds precisely to the complex structure on the irreducible components of a stable curve \emph{which contain no marked points}.

Let $C$ be a stable curve. We can use the Jenkins-Strebel theory to associate to $C$ a piece of combinatorial data which describes some of the complex structure of $C$. We start by deleting the nodal singularities of $C$, then we apply the Jenkins-Strebel theory to each connected component of the resulting surface in the same way as has already been discussed; that is to say that we choose a list of perimeters for the marked points and consider the critical graph of the corresponding Jenkins-Strebel differential on each connected component. The caveat here of course is that we cannot apply the Jenkins-Strebel theory to those components without any marked points. For this reason we collapse these components of the curve and label the resulting nodal singularity by the (arithmetic) genus of the collapsed component; we call this number the \emph{genus defect}.

This surface now has a connected graph embedded inside it defined by the union of the critical graphs of the Jenkins-Strebel differentials. Each vertex of this graph is labeled by the nonnegative integer corresponding to the genus of the collapsed component (or zero if no component was collapsed in forming that vertex). The edges issuing from each vertex are partitioned into groups called \emph{cycles}, the number of cycles being equal to the number of components meeting that vertex. The edges belonging to any particular cycle all lie on one of the components and therefore pick up a cyclic ordering as before. We call the graph resulting from this procedure a \emph{stable ribbon graph}. If we further decorate every edge of this graph with the positive real number defined by the length of this edge in the metric determined by the Jenkins-Strebel differential, we arrive at the definition of a \emph{stable metric ribbon graph}.

\begin{figure}[htp]
\centering
\includegraphics{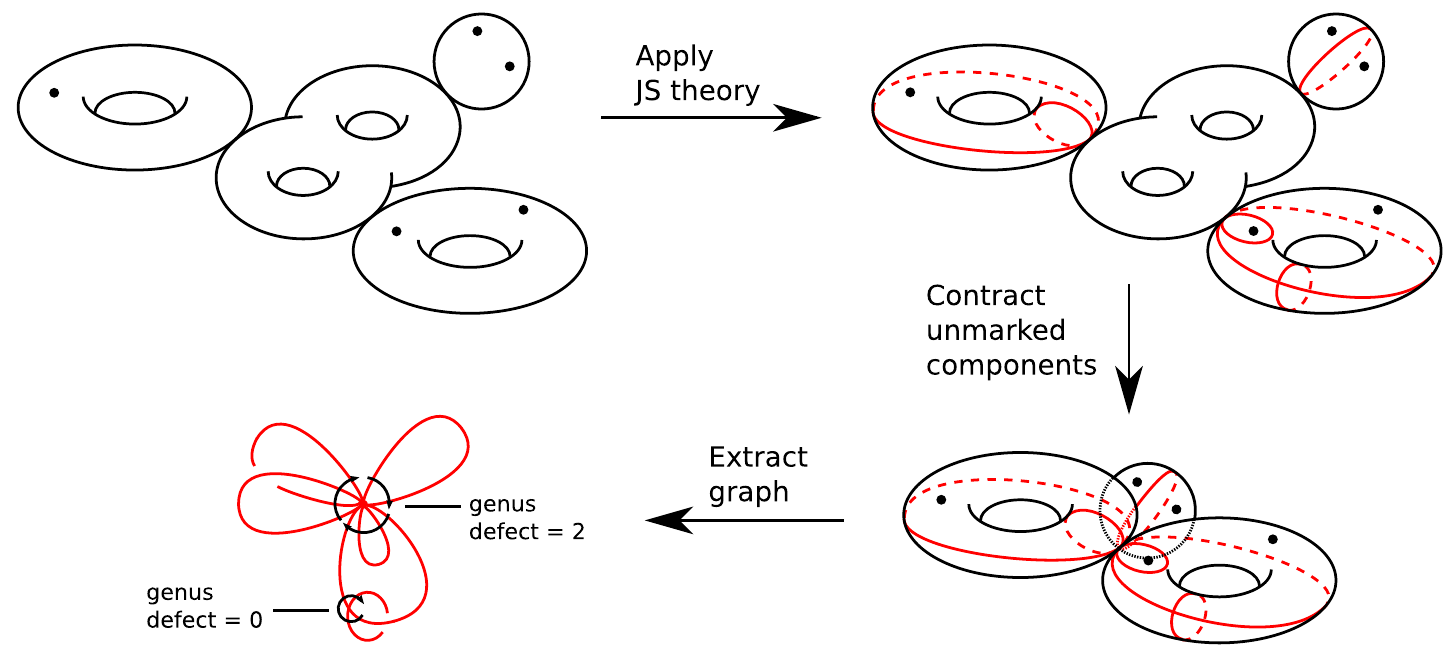}
\caption{Extracting a stable ribbon graph from a stable curve.}
\end{figure}

Hence we have seen that to any stable curve we can associate a certain piece of combinatorial information called a stable metric ribbon graph. Unfortunately, as has already been mentioned, it is not possible to reconstruct the original curve from this data since the Jenkins-Strebel theory cannot be used to recover the complex structure on the components without marked points. To this end we introduce the following equivalence relation on $\dmcmp$: consider two stable curves $C$ and $C'$, collapse those components without marked points and label the resulting nodal singularities with their arithmetic genus; then $C$ and $C'$ are equivalent if the resulting curves are biholomorphic through a mapping which preserves the genus defect parameters. This is the compactification of the moduli space defined by Kontsevich in \cite{kontairy}.

Let us denote the quotient of the Deligne-Mumford compactification by this equivalence relation by $\kcmp$. Using the above procedure it is clear that we can associate a stable metric ribbon graph to any point in $\kcmp\times\Delta^{\circ}_{n-1}$. Moreover, given any stable metric ribbon graph, a standard gluing construction allows us to reconstruct the point in $\kcmp\times\Delta^{\circ}_{n-1}$. This leads to an orbi-cellular decomposition of $\kcmp\times\Delta^{\circ}_{n-1}$; two points in $\kcmp\times\Delta^{\circ}_{n-1}$ belong to the same cell if and only if their corresponding stable ribbon graphs are isomorphic. We summarise this in the following theorem, due to Kontsevich \cite{kontairy}:

\begin{theorem} \label{thm_orbikmod}
The one-point compactification of the decorated moduli space $\kcmp\times\Delta^{\circ}_{n-1}$ is an orbi-cellular complex whose orbi-cells are indexed by stable ribbon graphs (and one 0-cell for the point). An orbi-cell $E$ lies on the boundary of another orbi-cell $E'$ if and only if the stable ribbon graph corresponding to $E$ is obtained from the stable ribbon graph corresponding to $E'$ by contracting some of the edges.
\end{theorem}
\noproof

\begin{rem}
Of course, one must describe exactly how the edges are contracted for a stable ribbon graph. This will be done in Section \ref{sec_ribgraph}.
\end{rem}

The application of the Jenkins-Strebel theory, as outlined above, can in fact be extended to give an orbi-cellular decomposition of a compactification of the decorated moduli space $\mspc\times\Delta_{n-1}$. This compactification slightly generalises Kontsevich's original construction and was introduced by Looijenga in \cite{looi}.

Let $C$ be a stable curve and let $p_1,\ldots,p_n$ be a list of perimeters for the marked points, some of which (but not all) are allowed to vanish. Again, we will use the Jenkins-Strebel theory to associate a type of graph to this piece of information. We begin by deleting the nodal singularities and puncturing the surface at those marked points with \emph{vanishing perimeters}. Next, we apply the Jenkins-Strebel theory to each connected component of this surface. Of course, the caveat here is that we cannot apply the Jenkins-Strebel theory to a component which has no marked points, or to a component for which all of the marked points have vanishing perimeters; therefore we collapse these components and label the corresponding nodal singularity by both the arithmetic genus of the collapsed component and the number of marked points that it contains. The first number is referred to as the \emph{genus defect}, as before, and the second number is referred to as the \emph{boundary defect}.

This surface has a connected graph lying embedded inside of it determined by the critical graphs of the Jenkins-Strebel differentials. The structure of this graph is exactly the same as before except that each vertex is decorated by an additional parameter coming from the boundary defect. For this reason, we will also refer to such a graph as a \emph{stable ribbon graph}.

\begin{figure}[htp]
\centering
\includegraphics{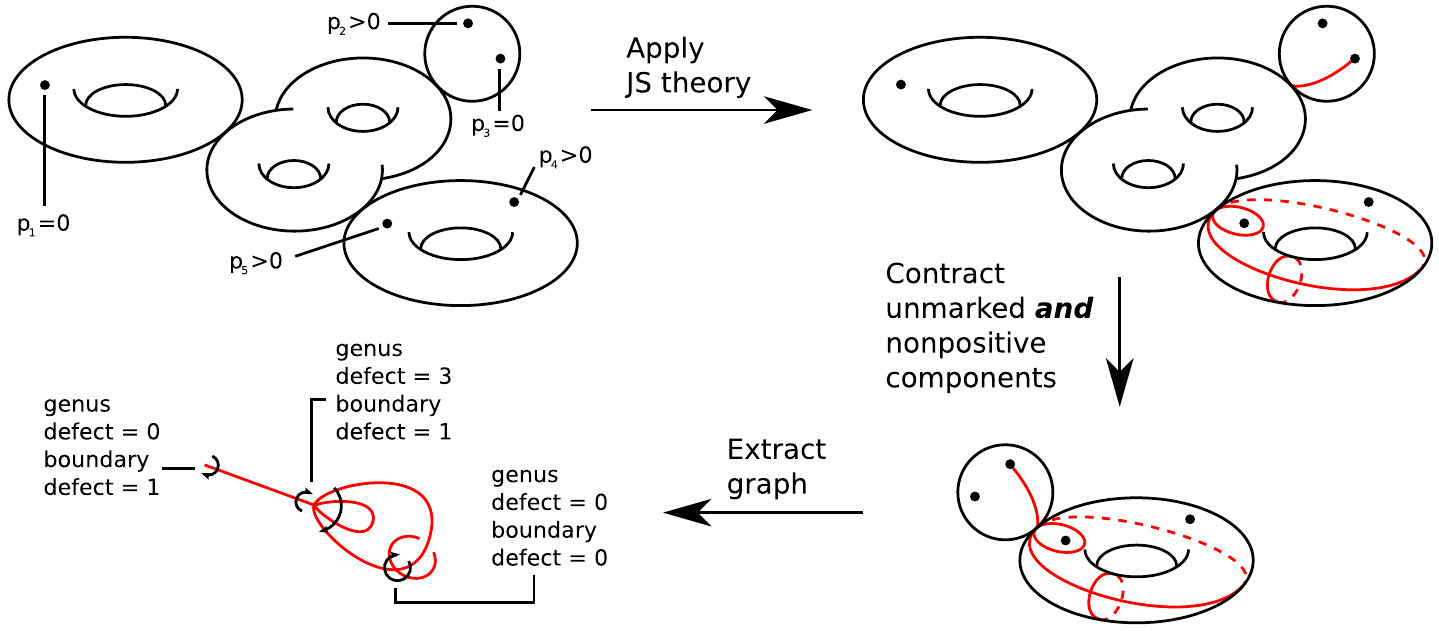}
\caption{Extracting a stable ribbon graph from a stable curve with \emph{vanishing perimeters}.}
\end{figure}

Now let us describe what kind of moduli space we have constructed an orbi-cellular decomposition of. Here we see that we are losing even more information than before; we are losing not just the complex structure of those components without marked points, but also the complex structure on those components whose marked points all have vanishing perimeters.

To this end we introduce the following equivalence relation on $\dmcmp\times\Delta_{n-1}$: two decorated curves $(C,p_1,\ldots,p_n)$ and $(C',p'_1,\ldots,p'_n)$ are equivalent if $p_i=p'_i$ for all $i$ and the curves given by collapsing both the components without marked points and those components whose marked points all have vanishing perimeters, are biholomorphic through a mapping which preserves both the genus and the boundary defect. We will denote the quotient of $\dmcmp\times\Delta_{n-1}$ by this equivalence relation by $\lcmp$. We have the following theorem, due to Looijenga \cite{looi}:

\begin{theorem} \label{thm_orbilmod}
The moduli space $\lcmp$ is an orbi-cellular complex whose orbi-cells are indexed by stable ribbon graphs. An orbi-cell $E$ lies on the boundary of another orbi-cell $E'$ if and only if the stable ribbon graph corresponding to $E$ is obtained from the stable ribbon graph corresponding to $E'$ by contracting some of the edges.
\end{theorem}
\noproof

\section{Noncommutative geometry and Lie bialgebras} \label{sec_ncgeom}

In this section we recall the basic framework of noncommutative symplectic geometry as defined by Kontsevich in \cite{kontsympgeom}. After recalling how the Lie algebra of noncommutative Hamiltonians on a symplectic vector space is defined, we introduce a definition for the divergence of a noncommutative vector field. We then use our definition for the divergence of a noncommutative vector field to give Kontsevich's Lie algebra of noncommutative Hamiltonians the extra structure of a \emph{Lie bialgebra}.

\subsection{Noncommutative differential geometry}

We begin by recalling the definition of noncommutative (polynomial) differential forms and the corresponding definition of the de Rham complex. The first place we must start is the definition of noncommutative 1-forms.

\begin{defi}
Let $V$ be a vector space. The module of noncommutative 1-forms $\Omega^1(V)$ is defined as
\[ \Omega^1(V):= T(V^*) \otimes T^+(V^*). \]
$\Omega^1(V)$ has the structure of a $T(V^*)$-bimodule via the actions
\begin{align*}
& a \cdot (x\otimes y):= ax \otimes y, \\
& (x \otimes y) \cdot a:= x \otimes ya - xy \otimes a;
\end{align*}
for $a,x \in T(V^*)$ and $y \in T^+(V^*)$.

Let $d:T(V^*) \to \Omega^1(V)$ be the map given by the formulae
\begin{displaymath}
\begin{array}{ll}
d(x):= 1 \otimes x, & x \in T^+(V); \\
d(x):= 0, & x \in \gf.
\end{array}
\end{displaymath}
The map $d$ thus defined is a derivation of degree zero.
\end{defi}

The noncommutative 1-forms are used to construct the algebra of all noncommutative forms as follows:

\begin{defi} \label{def_forms}
Let $V$ be a vector space and let $A:=T(V^*)$. The algebra of noncommutative forms $\forms{V}$ is defined as
\[ \forms{V}:=T_A\left[\Pi\Omega^1(V)\right]=A \oplus \bigoplus_{i=1}^\infty \underbrace{\Pi\Omega^1(V)\underset{A}{\otimes} \ldots \underset{A}{\otimes} \Pi\Omega^1(V)}_{i \text{ factors}}. \]

Since $\Omega^1(V)$ is an $A$-bimodule, $\forms{V}$ has the structure of an  associative algebra whose multiplication is the standard associative multiplication on the tensor algebra $T_A\left[\Pi\Omega^1(V)\right]$. The map $d:T(V^*) \to \Omega^1(V)$ lifts uniquely to a map $d:\forms{V} \to \forms{V}$ which gives $\forms{V}$ the structure of a differential graded algebra.
\end{defi}

It is possible to introduce analogues of the Lie derivative and contraction operator on the algebra of noncommutative forms, which is done as follows:

\begin{defi} \label{def_operate}
Let $V$ be a vector space and let $\xi:T(V^*) \to T(V^*)$ be a vector field:
\begin{enumerate}
\item
We can define a vector field $L_\xi: \forms{V} \to \forms{V}$, called the Lie derivative, by the formulae:
\begin{align*}
L_\xi(x)&:=\xi(x), \\
L_\xi(dx)&:=(-1)^{|\xi|}d(\xi(x)); \\
\end{align*}
for any $x \in T(V^*)$.
\item
We can define a vector field $i_\xi:\forms{V} \to \forms{V}$, called the contraction operator, by the formulae:
\begin{align*}
i_\xi(x)&:=0, \\
i_\xi(dx)&:=\xi(x); \\
\end{align*}
for any $x \in T(V^*)$.
\end{enumerate}
\end{defi}

It turns out that the algebra of noncommutative forms is not the right thing to consider in the framework of noncommutative geometry. Instead, we must consider its quotient by the submodule of commutators.

\begin{defi}
Let $V$ be a vector space. The de Rham complex $\drham{V}$ is defined as
\[ \drham{V}:=\frac{\forms{V}}{\left[\forms{V},\forms{V}\right]}. \]
The differential on $\drham{V}$ is induced by the differential on $\forms{V}$ defined in Definition \ref{def_forms} and is similarly denoted by $d$.
\end{defi}

The definition of the Lie derivative and contraction operator pass naturally to this quotient to give Lie and contraction operators on the de Rham complex $\drham{V}$. Of course, quotienting out by the submodule of commutators means that the de Rham complex is no longer an algebra. By an abuse of terminology, we will continue to refer to elements of $\drham{V}$ as differential forms.

\subsection{Noncommutative symplectic geometry}

Having now recalled the basic framework of noncommutative differential geometry, we proceed to introduce the relevant terminology for noncommutative \emph{symplectic} geometry, as outlined in \cite{ginzsympgeom}, \cite{kontsympgeom} and \cite{hamcohom}. We start with the definition of a symplectic form and a symplectic vector field.

\begin{defi} \label{def_sympform}
Let $V$ be a vector space and $\omega \in \drham[2]{V}$ be any 2-form. We say that $\omega$ is a \emph{symplectic form} if:
\begin{enumerate}
\item
it is a closed form, that is to say that $d\omega = 0$;
\item
it is nondegenerate, that is to say that the following map is bijective;
\begin{equation}
\begin{array}{ccc}
\Der[T(V^*)] & \to & \drham[1]{V}, \\
\xi & \mapsto & i_\xi(\omega).
\end{array}
\end{equation}
\end{enumerate}
\end{defi}

\begin{defi}
Let $V$ be a vector space and let $\omega \in \drham[2]{V}$ be a symplectic form. We say a vector field $\xi:T(V^*) \to T(V^*)$ is a \emph{symplectic vector field} if $L_\xi(\omega)=0$.
\end{defi}

In what follows we will only consider \emph{constant} symplectic forms. A constant 2-form is a 2-form $\omega\in\drham[2]{V}$ which can be written in the form 
\begin{equation} \label{eqn_constwoform}
\omega=\sum_i dx_i dy_i
\end{equation}
for some functions $x_i,y_i\in V^*$.

There is a one-to-one correspondence between constant 2-forms and skew-symmetric bilinear forms: given any constant 2-form $\omega$ as in \eqref{eqn_constwoform} we define the corresponding bilinear form $\innprod$ by the formula,
\begin{equation} \label{eqn_bilinform}
\langle a,b \rangle:=\sum_i (-1)^x_i[x_i(a)y_i(b)-(-1)^{ab}y_i(a)x_i(b)].
\end{equation}
Furthermore, the symplectic form $\omega$ is nondegenerate if and only if the bilinear form $\innprod$ is nondegenerate.

Any nondegenerate bilinear form $\innprod$ on a vector space yields a nondegenerate bilinear form $\innprod^{-1}$ on the dual space, defined by simply identifying the space with its dual. If we assume that
\[ \underbrace{x_1,\ldots,x_k}_{even};\underbrace{\xi_1,\ldots,\xi_k}_{odd}\in V^* \]
is a system of coordinates on $V$ and that our symplectic form $\omega$ is given by the formula
\[ \omega=\sum_{i=1}^k dx_i d\xi_i \]
then we have the following formula for $\innprod^{-1}$:
\[ \langle x_i,\xi_j \rangle^{-1} = \langle \xi_j,x_i \rangle^{-1} = \delta_{ij}. \]

\subsection{Lie algebras of noncommutative vector fields}

In this section we recall how to define a Lie algebra structure on the space of noncommutative 0-forms when the underlying manifold is equipped with a symplectic form. In what follows, we assume that our symplectic form is \emph{odd}, although an analogous treatment applies when our symplectic form is even.

Let $(V,\omega)$ be a symplectic vector space whose symplectic form is \emph{odd}. To any 0-form $a\in\drham[0]{V}$ we can associate a certain symplectic vector field $\alpha\in\Der[T(V^*)]$. This symplectic vector field is uniquely specified by the equation
\begin{equation} \label{eqn_formfield}
da=i_{\alpha}(\omega).
\end{equation}
This correspondence allows us to define an odd Lie bracket, also known as an anti-bracket, on the space $\mathfrak{h}[V]:=\drham[0]{V}$.

\begin{defi} \label{def_bracket}
Given a symplectic vector space $(V,\omega)$ as above, we define a bracket
\[ \{-,-\}:\mathfrak{h}[V]\otimes\mathfrak{h}[V]\to\mathfrak{h}[V] \]
of \emph{odd} degree by the formula
\[ \{a,b\}:=L_\alpha(b) \]
\end{defi}

\begin{prop}
The bracket $\{-,-\}$ on $\mathfrak{h}[V]$ is an odd Lie bracket, that is to say that the bracket
\[ [-,-]:\Pi\mathfrak{h}[V]\otimes\Pi\mathfrak{h}[V]\to\Pi\mathfrak{h}[V] \]
given by the formula $\Pi\circ [-,-] = \{-,-\}\circ(\Pi\otimes\Pi)$ is a Lie bracket.
\end{prop}
\noproof

The proof follows as a result of standard identities for the operators introduced in Definition \ref{def_operate}. The Lie algebra structure on $\mathfrak{h}[V]$ corresponds precisely under \eqref{eqn_formfield} to the usual commutator bracket of symplectic vector fields.

Let $a_1,\ldots,a_n;b_1,\ldots,b_m\in V^*$ be linear functions. An explicit formula for the Lie bracket $\{-,-\}$ is
\begin{equation} \label{eqn_bracket}
\{a_1\cdots a_n,b_1\cdots b_m\} = \sum_{i=1}^n\sum_{j=1}^m (-1)^p \langle a_i,b_j \rangle^{-1} (z_{n-1}^{i-1}\cdot [a_1\cdots \hat{a_i} \cdots a_n]) (z_{m-1}^{j-1}\cdot [b_1\cdots \hat{b_j} \cdots b_m]),
\end{equation}
where $z_k$ denotes the permutation $(k \, k-1\ldots 2 \, 1)$ and
\[ p:=|a_i|(|a_1|+\ldots+|a_{i-1}|) + |b_j|(|a_1|+\ldots+|a_n|+|b_1|+\ldots+|b_{j-1}|). \]

\subsection{A Lie bialgebra structure on the space of noncommutative 0-forms}

In this section we give a definition for the divergence of a noncommutative vector field and use this to construct a Lie cobracket on the space of noncommutative 0-forms on a symplectic vector space. We show that combining this structure with the Lie algebra structure described in the last section gives the space of 0-forms the structure of an involutive Lie bialgebra.

Our definition for the divergence of a noncommutative vector field has a slightly curious appearance as it lands in the second tensor power of 0-forms, instead of landing in 0-forms as is usual in commutative geometry. Nevertheless, it does lift the ordinary definition for the divergence of a commutative vector field and satisfies an important identity which is the analogue of the classical formula for the divergence of a commutator of two vector fields.

\begin{defi}
Let $V$ be a finite-dimensional vector space with coordinates $x_1,\ldots,x_n$ and let $f_1,\ldots,f_k\in V^*$ be linear functions. We define the divergence of the vector field $\xi:=(f_1\cdots f_k)\partial_{x_i}$ by the formula
\begin{equation} \label{eqn_divergence}
\nabla(\xi):=\sum_{i=1}^k (-1)^{x_i(f_i+\ldots+f_k)} \partial_{x_i}(f_i)\cdot[(f_1\cdots f_{i-1})\otimes (f_{i+1}\cdots f_k)].
\end{equation}
By extending \eqref{eqn_divergence} linearly, we arrive at the definition for the divergence
\[\nabla:\Der[T(V^*)]\to\drham[0]{V}\otimes\drham[0]{V}.\]
\end{defi}

It is easy to see that the above definition of divergence is independent of the choice of coordinates. It satisfies the following formula for the commutator of two vector fields:

\begin{lemma} \label{lem_divcom}
Let $V$ be a finite-dimensional vector space and $\xi,\gamma\in\Der[T(V^*)]$ be vector fields, then
\[ \nabla([\xi,\gamma]) = (L_{\xi}\otimes 1 + 1\otimes L_{\xi})[\nabla(\gamma)] - (-1)^{\xi\gamma} (L_{\gamma}\otimes 1 + 1\otimes L_{\gamma})[\nabla(\xi)].\]
\end{lemma}

\begin{proof}
The proof follows by direct calculation.
\end{proof}

In particular, this formula guarantees that the subspace of noncommutative vector fields with vanishing divergence forms a Lie subalgebra of $\Der[T(V^*)]$. It also corresponds to the compatibility condition between the Lie bracket on 0-forms and the Lie cobracket on 0-forms, which we will define next. First of all, however, let us recall the definition of an involutive Lie bialgebra.

\begin{defi}
A Lie bialgebra is a vector space $\mathfrak{g}$ together with the structures of a Lie bracket
\[ [-,-]:\mathfrak{g}\otimes\mathfrak{g}\to\mathfrak{g} \]
and a Lie cobracket
\[ \Delta:\mathfrak{g}\to\mathfrak{g}\otimes\mathfrak{g}\]
such that the following compatibility condition is satisfied:
\begin{equation} \label{eqn_cocycle}
\Delta([x,y])=[x,\Delta(y)]-(-1)^{xy}[y,\Delta(x)].
\end{equation}
Furthermore, we say that $\mathfrak{g}$ is an \emph{involutive} Lie bialgebra if the following additional condition is satisfied:
\begin{equation} \label{eqn_involutive}
[-,-]\circ\Delta=0.
\end{equation}
\end{defi}

\begin{rem}
Note that for $x,y,z\in\mathfrak{g}$ we define
\[ [x,y\otimes z]:=[x,y]\otimes z + (-1)^{xy}y\otimes [x,z]. \]
\end{rem}

Now we use our definition for the divergence of a noncommutative vector field to define an involutive Lie bialgebra structure on the space of 0-forms on a symplectic vector space.

\begin{defi}
Let $(V,\omega)$ be a symplectic vector space whose symplectic form is odd. We define a diagonal $\Delta$ on $\mathfrak{h}[V]:=\drham[0]{V}$ of \emph{odd} degree by the following commutative diagram:
\[\xymatrix{ \mathfrak{h}[V] \ar[rr]^-{\Delta} \ar[rd]^{a\mapsto\alpha} && \mathfrak{h}[V]\otimes\mathfrak{h}[V] \\ & \Der[T(V^*)] \ar[ru]^{\frac{1}{2}\nabla} }\]
where the map in the lower left corner is that defined by Equation \eqref{eqn_formfield}.
\end{defi}

\begin{prop} \label{prop_bialgebra}
The diagonal $\Delta$ on $\mathfrak{h}[V]$, together with the bracket $\{-,-\}$ on $\mathfrak{h}[V]$ described in Definition \ref{def_bracket} give $\mathfrak{h}[V]$ (or, more precisely, its parity reversion $\Pi\mathfrak{h}[V]$) the structure of an involutive Lie bialgebra.
\end{prop}

\begin{proof}
An explicit formula for the cobracket $\Delta$ is given by the following: let $a_1,\ldots,a_n\in V^*$ be linear functions, then
\begin{equation} \label{eqn_cobracket}
\Delta(a_1\cdots a_n) = \frac{1}{2}\sum_{i<j}(-1)^p\langle a_i,a_j \rangle^{-1}[1+(1 \, 2)]\cdot[ (a_{i+1}\cdots a_{j-1}) \otimes (a_{j+1}\cdots a_n a_1\cdots a_{i-1})];
\end{equation}
where
\begin{displaymath}
\begin{split}
p:= & |a_i|(|a_1|+\ldots+|a_i|) + |a_j|(|a_1|+\ldots+|a_j|) \\
& + (|a_1|+\ldots+|a_{i-1}|)(|a_{i+1}|+\ldots+|a_{j-1}|+|a_{j+1}|+\ldots+|a_n|).
\end{split}
\end{displaymath}
Using this formula, one can verify directly both the coJacobi identity and the involutivity condition. The compatibility condition between the bracket and the cobracket follows as a direct consequence of Lemma \ref{lem_divcom}.
\end{proof}

\section{A two-parameter family of differential graded Lie algebras} \label{sec_dgla}

In this section we will use the framework of noncommutative geometry defined in the previous section to construct a two-parameter family of differential graded Lie algebras. This two-parameter family will be the central object of the main theorem formulated in Section \ref{sec_main}.

\subsection{The Chevalley-Eilenberg complex}

In this section we recall the definition of the Chevalley-Eilenberg complex of a differential graded Lie algebra as well as some of its basic properties. In this paper, we will only be interested in Lie algebras equipped with an antibracket; that is to say that we assume that the Lie bracket is a map of \emph{odd} degree.

\begin{defi}
Let $\mathfrak{g}$ be a differential graded Lie algebra whose bracket is an \emph{odd} map. The Chevalley-Eilenberg complex of $\mathfrak{g}$, denoted by $\ce{\mathfrak{g}}$, is the complex whose underlying vector space is the symmetric algebra on $\mathfrak{g}$:
\[ \ce{\mathfrak{g}}:=S(\mathfrak{g})=\bigoplus_{n=0}^{\infty} (\mathfrak{g}^{\otimes n})_{S_n}. \]
The differential $\delta:\ce{\mathfrak{g}}\to\ce{\mathfrak{g}}$ is defined by the formula,
\begin{displaymath}
\begin{split}
\delta(g_1\cdots g_n):=&\sum_{1\leq i<j\leq n} (-1)^p\{g_i,g_j\}\cdot g_1\cdots \hat{g_i} \cdots \hat{g_j} \cdots g_n \\
& + \sum_{1\leq i\leq n} (-1)^q d(g_i)\cdot g_1\cdots \hat{g_i} \cdots  g_n;
\end{split}
\end{displaymath}
where
\begin{displaymath}
\begin{split}
p:= & |g_i|(|g_1|+\ldots+|g_{i-1}|) + |g_j|(|g_1|+\ldots+|g_{j-1}|) + |g_i||g_j|, \\
q:= & |g_i|(|g_1|+\ldots+|g_{i-1}|)
\end{split}
\end{displaymath}
and $d$ is the differential on $\mathfrak{g}$. The homology of this complex is known as the Chevalley-Eilenberg homology of the differential graded Lie algebra $\mathfrak{g}$ and is denoted by $\hce{\mathfrak{g}}$.
\end{defi}

In fact the Chevalley-Eilenberg complex has much more structure than simply that of a complex. It has a commutative multiplication
\[ - \cdot -:\ce{\mathfrak{g}}\otimes\ce{\mathfrak{g}}\to\ce{\mathfrak{g}} \]
coming from the canonical multiplication in the symmetric algebra and we can also equip it with an \emph{odd} bracket
\[ \{-,-\}:\ce{\mathfrak{g}}\otimes\ce{\mathfrak{g}}\to\ce{\mathfrak{g}} \]
by simply extending the bracket on $\mathfrak{g}$ according to the Leibniz rule. This bracket is given by the formula
\begin{equation} \label{eqn_bvbracket}
\{g_1\cdots g_n,h_1\cdots h_m\}=\sum_{i=1}^n\sum_{j=1}^m (-1)^p\{g_i,h_j\}\cdot g_1\cdots \hat{g_i} \cdots g_n\cdot h_1\cdots \hat{h_j} \cdots h_m,
\end{equation}
where
\[ p:=|g_i|(|g_1+\ldots+|g_{i-1}|) + |h_j|(|g_1|+\ldots+|g_n|+|h_1|+\ldots+|h_{j-1}|) + |g_i||h_j|. \]
The appropriate terminology for this type of algebraic structure is a \emph{Batalin-Vilkovisky algebra}, whose definition we will now recall.

\begin{defi}
A Batalin-Vilkovisky algebra is a vector space $W$ equipped with:
\begin{enumerate}
\item
a differential $d:W\to W$,
\item
a commutative product $- \cdot -:W\otimes W\to W$ of \emph{even} degree and
\item
a Lie bracket $\{-,-\}:W\otimes W\to W$ of \emph{odd} degree.
\end{enumerate}

These structures must satisfy the following axioms:
\begin{enumerate}
\item
The bracket and product must satisfy the Leibniz rule; that is to say that for all $a,b,c\in W$,
\[ \{a,b\cdot c\} = \{a,b\}\cdot c + (-1)^{(a+1)b}b\cdot\{a,c\}. \]
\item
The differential should be a derivation of the Lie bracket; that is to say that for all $a,b\in W$,
\[ d(\{a,b\})+\{d(a),b\}+(-1)^a\{a,d(b)\}=0.\]
\item
For all $a,b\in W$,
\[d(a\cdot b)=d(a)\cdot b + (-1)^a a\cdot d(b) + \{a,b\}.\]
\end{enumerate}
\end{defi}

In fact, the second axiom is a consequence of the third axiom. It is a standard fact, which can be verified directly, that the Chevalley-Eilenberg complex with the algebraic structures defined above is a Batalin-Vilkovisky algebra; in particular, it is a differential graded Lie algebra.

In Section \ref{sec_main}, we will need to consider a minor variant of Chevalley-Eilenberg homology known as \emph{relative} Chevalley-Eilenberg homology. We now recall its definition.

\begin{defi}
Let $\mathfrak{g}$ be any differential graded Lie algebra and let $\mathfrak{h}\subset\mathfrak{g}$ be an arbitrary differential graded Lie subalgebra. Note that $\mathfrak{h}$ acts on $S(\mathfrak{g}/\mathfrak{h})$ canonically, as $\mathfrak{g}/\mathfrak{h}$ is an $\mathfrak{h}$-module. The relative Chevalley-Eilenberg complex $\rce{\mathfrak{g}}{\mathfrak{h}}$ is given by taking the coinvariants of this action:
\[ \rce{\mathfrak{g}}{\mathfrak{h}}:=S(\mathfrak{g}/\mathfrak{h})_{\mathfrak{h}}. \]
One can check that the differential $\delta$ on $\ce{\mathfrak{g}}$ induces a well-defined differential on $\rce{\mathfrak{g}}{\mathfrak{h}}$, also denoted by $\delta$. The homology of this complex is called the \emph{relative} Chevalley-Eilenberg homology of $\mathfrak{g}$ modulo $\mathfrak{h}$ and is denoted by $\hrce{\mathfrak{g}}{\mathfrak{h}}$.
\end{defi}

\subsection{Construction of the two-parameter family} \label{sec_constructdgla}

In this section we construct the two-parameter family of differential graded Lie algebras that will play the central role in the main theorem of Section \ref{sec_main}. Here we exploit a standard construction which produces a differential graded Lie algebra from any involutive Lie bialgebra. We apply this construction to the involutive Lie bialgebra structure we defined in Section \ref{sec_ncgeom} on the space of 0-forms on a symplectic vector space.

We begin by recalling the details of this construction. Let $(V,\omega)$ be a symplectic vector space whose symplectic form is \emph{odd} and let $\mathfrak{h}:=\mathfrak{h}[V]$ be the involutive Lie bialgebra of Proposition \ref{prop_bialgebra}. Consider the Chevalley-Eilinberg complex $(\ce{\mathfrak{h}},\delta)$. We know from the results of the previous section that $\ce{\mathfrak{h}}$ is a differential graded Lie algebra. There is a way to include the Lie bialgebra structure on $\mathfrak{h}$ in the Chevalley-Eilenberg complex $\ce{\mathfrak{h}}$ such that the Chevalley-Eilenberg complex retains the structure of a differential graded Lie algebra. We can define a map
\[ \Delta:\ce{\mathfrak{h}}\to\ce{\mathfrak{h}} \]
from the Lie cobracket $\Delta:\mathfrak{h}\to\mathfrak{h}\otimes\mathfrak{h}$ by simply extending the cobracket using the Leibniz rule:
\[\Delta(h_1\cdots h_n):=\sum_{i=1}^n (-1)^p\Delta(h_i)\cdot h_1\cdots \hat{h_i} \cdots h_n,\]
where $p:=|h_i|(|h_1|+\ldots+|h_{i-1}|)$.

Now we tensor the Chevalley-Eilenberg complex $\ce{\mathfrak{h}}$ with the free polynomial algebra in one generator $\gamma$ and equip it with a deformed differential:

\begin{lemma}
For any symplectic vector space $(V,\omega)$, the tensor product of the Chevalley-Eilenberg complex of $\mathfrak{h}:=\mathfrak{h}[V]$ with the free polynomial algebra in one variable $\gamma$
\[ \mathfrak{l}:=\gf[\gamma]\otimes\ce{\mathfrak{h}}\]
is a differential graded Lie algebra when we equip it with the differential $d:=\gamma\cdot\delta + \Delta$.
\end{lemma}

\begin{proof}
The Lie bracket on the Chevalley-Eilenberg complex of $\mathfrak{h}$ extends naturally to $\mathfrak{l}$ as a trivial deformation. We already know that $\delta$ is a derivation of this Lie bracket; that $\Delta$ is also a derivation follows from the compatibility condition \eqref{eqn_cocycle} between the bracket and the cobracket on $\mathfrak{h}$.

Next we show that $d^2=0$. Since we already know that $\delta^2=0$ and since the condition $\Delta^2=0$ follows from the coJacobi identity, this is equivalent to $[\delta,\Delta]=0$. The compatibility condition \eqref{eqn_cocycle} between the bracket and cobracket ensures that $[\delta,\Delta]$ is a derivation of the commutative product and the involutivity constraint \eqref{eqn_involutive} guarantees that it is zero on the generators and hence zero everywhere.
\end{proof}

Recall that the underlying vector space of $\mathfrak{h}$ is
\[ \drham[0]{V}:=\bigoplus_{i=0}^\infty ([V^*]^{\otimes i})_{\mathbb{Z}/i\mathbb{Z}}. \]
Let us introduce the notation $\mathfrak{h}_{\geq n}$ for the subspace of 0-forms of order $\geq n$,
\[\mathfrak{h}_{\geq n}:=\bigoplus_{i=n}^\infty ([V^*]^{\otimes i})_{\mathbb{Z}/i\mathbb{Z}}.\]
With this notation the Lie algebra $\mathfrak{h}$ splits as a vector space,
\[\mathfrak{h}=\gf\oplus\mathfrak{h}_{\geq 1};\]
hence we see that the differential graded Lie algebra $\mathfrak{l}$ is really a two-parameter deformation
\[\mathfrak{l}=\gf[\gamma]\otimes S(\mathfrak{h})=\gf[\gamma,\nu]\otimes S(\mathfrak{h}_{\geq 1})\]
by identifying the symmetric algebra on the field $\gf$ with the free polynomial algebra in one generator $\nu$. One can check that the differential and the Lie bracket are actually $\gf[\gamma,\nu]$-linear, which follows as a simple consequence of the definitions.

Now for our purposes, this differential graded Lie algebra is not exactly what we want, hence we must embark on a technical description as to how it is to be altered. We need to modify it by cutting out some of the low order terms. Note that this issue also arises in Kontsevich's original paper \cite{kontsympgeom}, although there it is technically much more straightforward to deal with. The main problem is that the vertices of a stable ribbon graph which consist of only one cycle and have vanishing genus and boundary defect \emph{must be at least trivalent}.

We begin by noting that $\mathfrak{l}$ splits as a vector space
\[\mathfrak{l} = \gf[\gamma,\nu] \oplus \left(\gf[\gamma,\nu]\otimes\left[\bigoplus_{i=1}^\infty \left[(\mathfrak{h}_{\geq 1})^{\otimes i}\right]_{S_i}\right]\right). \]
Observing that the left-hand summand $\gf[\gamma,\nu]$ is a trivial ideal in our differential graded Lie algebra, we see that
\begin{equation} \label{eqn_dglacut}
\mathfrak{l}':=\gf[\gamma,\nu]\otimes\left[\bigoplus_{i=1}^\infty \left[(\mathfrak{h}_{\geq 1})^{\otimes i}\right]_{S_i}\right]
\end{equation}
inherits the structure of a differential graded Lie algebra when we quotient out $\mathfrak{l}$ by this ideal.

The differential graded Lie algebra that we are looking for sits inside \eqref{eqn_dglacut} as a differential graded Lie subalgebra. Note again that $\mathfrak{h}_{\geq 1}$ splits as a vector space
\[ \mathfrak{h}_{\geq 1} = V^* \oplus \mathfrak{h}_{\geq 2}. \]
Let $\gf_+[\gamma,\nu]$ denote the ideal of $\gf[\gamma,\nu]$ consisting of polynomials which \emph{vanish at the origin}, so that $\gf[\gamma,\nu]$ splits as
\[ \gf[\gamma,\nu]=\gf\oplus\gf_+[\gamma,\nu]. \]
Hence, \eqref{eqn_dglacut} splits as a sum of vector spaces.
\[ \mathfrak{l}' = V^*\oplus(\gf_+[\gamma,\nu]\otimes V^*)\oplus(\gf[\gamma,\nu]\otimes\mathfrak{h}_{\geq 2})\oplus\left(\gf[\gamma,\nu]\otimes\left[\bigoplus_{i=2}^\infty \left[(\mathfrak{h}_{\geq 1})^{\otimes i}\right]_{S_i}\right]\right) \]

If we throw away the leftmost summand $V^*$ then we arrive at the desired definition for our differential graded Lie algebra:

\begin{theorem} \label{thm_dgla}
Let $(V,\omega)$ be a symplectic vector space whose symplectic form is odd, then
\[ \Lambda_{\gamma,\nu}[V]:=(\gf_+[\gamma,\nu]\otimes V^*)\oplus(\gf[\gamma,\nu]\otimes\mathfrak{h}_{\geq 2})\oplus\left(\gf[\gamma,\nu]\otimes\left[\bigoplus_{i=2}^\infty \left[(\mathfrak{h}_{\geq 1})^{\otimes i}\right]_{S_i}\right]\right). \]
is a differential graded Lie algebra, whose differential is the unique differential induced by the deformed differential $d:=\gamma\cdot\delta + \Delta$.
\end{theorem}

\begin{proof}
It is simple to check that $\Lambda_{\gamma,\nu}[V]$ is a differential graded Lie subalgebra of \eqref{eqn_dglacut}.
\end{proof}

\begin{rem} \label{rem_linfields}
Now suppose that the vector space $V$ has dimension $n|n$ and consider the Lie algebra $\mathfrak{h}_{\geq 2}$. The subspace $S^2(V^*)\subset\mathfrak{h}_{\geq 2}$ of strictly quadratic Hamiltonians forms a Lie subalgebra which can be identified with the Lie algebra $\mathfrak{pe}[\mathbb{Q}^{n|n}]$ of linear endomorphisms of $V$ which preserve the bilinear form associated to the symplectic form $\omega$ by \eqref{eqn_bilinform}. This Lie algebra also sits inside $\Lambda_{\gamma,\nu}$ as a differential graded Lie subalgebra with trivial differential, simply by choosing the inclusion $S^2(V^*)\subset\mathfrak{h}_{\geq 2}\subset\Lambda_{\gamma,\nu}$ corresponding to the summand $\gf\subset\gf[\gamma,\nu]$ (note that $\mathfrak{h}_{\geq 2}$ is not itself a Lie subalgebra of $\Lambda_{\gamma,\nu}$, only $S^2()$). The Lie subalgebra $\mathfrak{pe}[\mathbb{Q}^{n|n}]$ and its invariant theory will play an important role in the proof of the main theorem in Section \ref{sec_main}.
\end{rem}

\subsection{A diagram of differential graded Lie algebras}

Since the differential graded Lie algebra that we defined in Theroem \ref{thm_dgla} is a two-parameter deformation, we may consider the family of differential graded Lie algebras defined by choosing specific values in $\gf$ for these deformation parameters; in particular, we may consider the differential graded Lie algebras defined by setting one or both of these deformation parameters to be zero.

More precisely, given a symplectic vector space $(V,\omega)$, there exists the unique structure of a differential graded Lie algebra on the spaces
\begin{displaymath}
\begin{split}
\Lambda_{\gamma}[V]:= & (\gf_+[\gamma]\otimes V^*)\oplus(\gf[\gamma]\otimes\mathfrak{h}_{\geq 2})\oplus\left(\gf[\gamma]\otimes\left[\bigoplus_{i=2}^\infty \left[(\mathfrak{h}_{\geq 1})^{\otimes i}\right]_{S_i}\right]\right), \\
\Lambda[V]:= & \mathfrak{h}_{\geq 2}\oplus\left( \bigoplus_{i=2}^\infty \left[(\mathfrak{h}_{\geq 1})^{\otimes i}\right]_{S_i}\right); \\
\end{split}
\end{displaymath}
such that the maps in the diagram
\[ \xymatrix{\Lambda_{\gamma,\nu}[V] \ar[r]^{\nu =0} & \Lambda_{\gamma}[V] \ar[r]^{\gamma =0} & \Lambda[V]} \]
given by setting the deformation parameters to zero are morphisms of differential graded Lie algebras.

Now consider the canonical projection
\[ \pi:\Lambda[V]\to\mathfrak{h}_{\geq 2}[V] \]
from $\Lambda[V]$ to the Lie algebra $\mathfrak{h}_{\geq 2}[V]$ of Hamiltonians of quadratic and higher order. One can check that this projection is a morphism of differential graded Lie algebras, where $\mathfrak{h}_{\geq 2}[V]$ is equipped with the trivial differential. Let us introduce the notation $\mathfrak{g}[V]$ for the Lie algebra $\mathfrak{h}_{\geq 2}[V]$. Combining the map $\pi$ with the morphism that was defined above by setting the deformation parameter $\gamma$ to equal zero, we arrive at the following important diagram of differential graded Lie algebras:
\begin{equation} \label{eqn_seqdglas}
\xymatrix{ \Lambda_{\gamma,\nu}[V] \ar[r] & \Lambda_{\gamma}[V] \ar[r] & \mathfrak{g}[V]. }
\end{equation}
Recall from Remark \ref{rem_linfields} that the Lie algebra $\mathfrak{pe}[\mathbb{Q}^{n|n}]$ sits inside each one of these differential graded Lie algebras as a subalgebra. A simple check reveals that the above diagram respects these embeddings of $\mathfrak{pe}[\mathbb{Q}^{n|n}]$. The geometric interpretation of diagram \eqref{eqn_seqdglas} in terms of moduli spaces of Riemann surfaces will be explained in Section \ref{sec_main}.

\section{The stable ribbon graph complex} \label{sec_ribgraph}

In this section we give a precise definition of the complex of stable ribbon graphs. The main point is to provide a precise description of what happens to the graph of nonclosed horizontal trajectories of a Jenkins-Strebel differential when the length of one of its edges tends to zero in the corresponding metric. The presentation provided here is mainly based on the accounts given by \cite{mondello} and \cite{zvonkine}.

\subsection{Stable ribbon graphs}

We begin with the formal definition of an (oriented) \emph{stable ribbon graph}.

\begin{defi} \label{def_stablegraph}
A stable ribbon graph is a set $\Gamma$ called the set of \emph{half-edges} together with the folllowing data:
\begin{enumerate}
\item
A partition of $\Gamma$ into pairs, denoted by $E(\Gamma)$, called the set of \emph{edges} of $\Gamma$.
\item
A partition of $\Gamma$, denoted by $V(\Gamma)$, called the set of \emph{vertices} of $\Gamma$. We will refer to the cardinality of a vertex $v\in V(\Gamma)$ as the \emph{valency} of $v$.
\item
For every vertex $v\in V(\Gamma)$, a further partition $C(v)$ of $v$ called the \emph{cycles} of $v$. Furthermore, we require that every cycle $c\in C(v)$ is endowed with a cyclic ordering of its elements.
\item
For every vertex $v\in V(\Gamma)$, a pair of nonnegative integers $g(v)$ and $n(v)$, called the \emph{genus defect} and \emph{boundary defect} respectively. Furthermore, if both $g(v)$ and $n(v)$ are equal to zero and $C(v)$ consists of a single cycle, then we impose the additional requirement that the vertex $v$ be at least trivalent.
\item
An ordering of the edges of $\Gamma$ modulo the action of the group of even permutations of the edges. This part of the data is called the \emph{orientation} on $\Gamma$.
\end{enumerate}
\end{defi}

\begin{rem}
Note that if we consider only those stable ribbon graphs for which every vertex has just a single cycle and for which both the genus and boundary defect are equal to zero, then we recover the usual definition of a \emph{ribbon graph}.
\end{rem}

There is a fairly obvious notion of isomorphism for stable ribbon graphs. Two stable ribbon graphs are isomorphic if there is a bijective mapping between their set of half-edges preserving the structures defined by items (1)-(5) of Definition \ref{def_stablegraph}.

Given a stable ribbon graph $\Gamma$, we can associate permutations $\sigma_0,\sigma_1,\sigma_\infty:\Gamma\to\Gamma$ defined as follows:
\begin{enumerate}
\item
$\sigma_1$ is defined as the fixed point free involution whose 2-cycles are the edges of $\Gamma$,
\item
$\sigma_0$ is defined as the permutation whose cycles are the cycles of the graph $\Gamma$,
\item
$\sigma_\infty:=\sigma_0^{-1}\sigma_1$.
\end{enumerate}

The cycles of the permutation $\sigma_\infty$ are called the \emph{perimeters} of the stable ribbon graph. This is because they trace out the set of nonclosed horizontal trajectories surrounding a given marked point on a Riemann surface; hence, the perimeters of a stable ribbon graph are in one-to-one correspondence with the set of marked points on the corresponding Riemann surface which have nonvanishing perimeters.

We need to pay special attention to those perimeters which are constituted of a single edge or loop. This is because the contraction of such an edge or loop corresponds to shrinking the length of this perimeter to zero, for which there are special combinatorial rules. A set of representative examples is provided by the following figures.
\begin{figure}[htp] \label{fig_boundaryloops}
\centering
\includegraphics{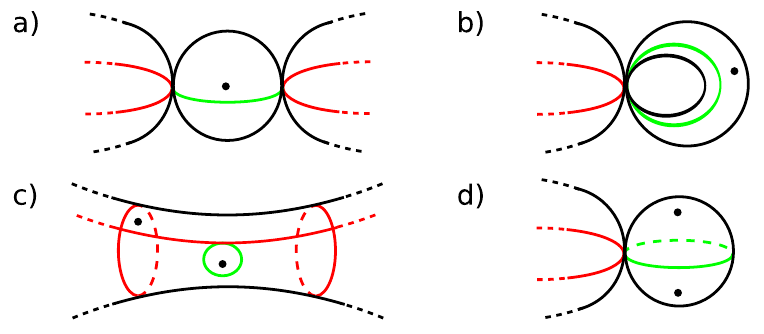}
\caption{Edges and loops completely surrounding a marked point.}
\end{figure}

We now describe the various combinatorial rules for contracting edges in a stable ribbon graph. They describe how the graph of nonclosed horizontal trajectories of a Jenkins-Strebel differential changes as we shrink the length of one of its edges.

\begin{defi} \label{def_contract}
Let $\Gamma$ be an oriented stable ribbon graph and let $e\in E(\Gamma)$ be an edge. We define the graph $\Gamma/e$ to be the graph obtained by contracting this edge according to the following rules:
\begin{enumerate}
\item \label{item_contractedge}
Suppose that $e$ is not a loop, in which case it joins distinct vertices $v_1,v_2\in V(\Gamma)$. These vertices are partitioned into cycles, so the endpoints of $e$ lie in distinct cycles $c_1\subset v_1$ and $c_2\subset v_2$. When the length of the edge $e$ shrinks to zero, the vertices $v_1$ and $v_2$ become joined, and the cycles $c_1$ and $c_2$ coalesce to form a new cycle with a naturally defined cyclic ordering. The genus and boundary defects for the vertices $v_1$ and $v_2$ are added to give the defects for the new vertex made from joining $v_1$ and $v_2$. The orientation is defined in an obvious way. All the other combinatorial structures elsewhere on the graph are left alone.
\begin{figure}[htp]
\centering
\includegraphics{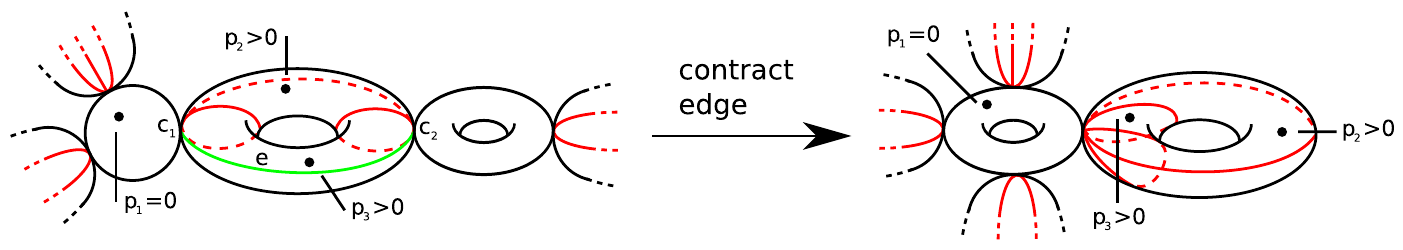}
\caption{Contracting an edge in a stable ribbon graph.}
\end{figure}

Note that when both $c_1$ and $c_2$ each consist of a single half-edge (cf. Figure \ref{fig_boundaryloops}(a)), $c_1$ and $c_2$ do not coalesce, but instead vanish and the boundary defect at the new vertex is defined to be the sum of the boundary defects of $v_1$ and $v_2$ \emph{plus one}. If, furthermore, $c_1$ and $c_2$ are the only cycles of $v_1$ and $v_2$, then the edge $e$ cannot be contracted.

\item \label{item_contractloop1}
Now suppose that $e$ is a loop, in which case both its endpoints lie in a single vertex $v$. Suppose furthermore, that they join distinct cycles $c_1,c_2\subset v$. As the length of the loop $e$ tends to zero, these cycles coalesce to form a single cycle as before. In so doing, a nonseparating double-point is formed on the topological surface corresponding to the vertex $v$, hence the genus defect of $v$ increases by one. No other combinatorial structures are changed.
\begin{figure}[htp]
\centering
\includegraphics{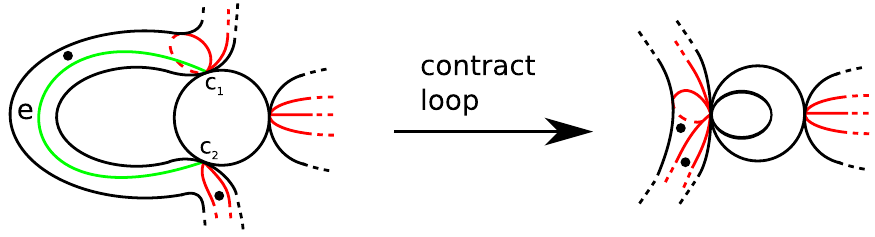}
\caption{Contracting a loop joining two distinct cycles.}
\end{figure}

As before, care must be taken when both $c_1$ and $c_2$ consist of a single half-edge (cf. Figure \ref{fig_boundaryloops}(b)). In this case $c_1$ and $c_2$ are annihilated and \emph{both} the genus \emph{and} the boundary defect are increased by one.

\item \label{item_contractloop2}
Finally, suppose that $e$ is again a loop, but that now both of its endpoints lie in the same cycle $c$ contained in some vertex $v$. Shrinking the length of this loop pinches the surface and a double-point is formed. The cycle $c$ splits up into two cycles $c_1$ and $c_2$, with naturally defined cyclic orderings. All the other combinatorial structures remain unchanged.
\begin{figure}[htp]
\centering
\includegraphics{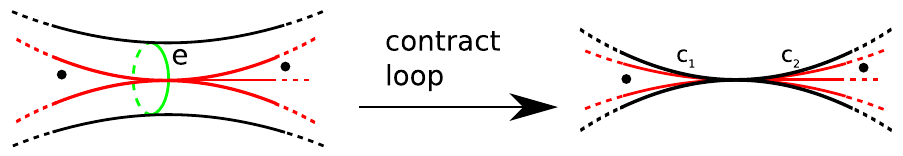}
\caption{Contracting a loop which joins a cycle to itself.}
\end{figure}

Again, care must be taken with this definition when the endpoints of $e$ lie next to each other in the cyclic ordering (cf. Figure \ref{fig_boundaryloops}(c)). In this case, the cycle $c$ does not split up, but the boundary defect is increased by one. Furthermore, if the cycle $c$ consists of just the two half-edges of $e$ (cf. Figure \ref{fig_boundaryloops}(d)), then the cycle $c$ is annihilated and the boundary defect actually increases by \emph{two}. Finally, if in the lattermost situation the vertex has no other cycles than $c$, then the loop $e$ cannot actually be contracted at all.
\end{enumerate}
\end{defi}

\subsection{The stable ribbon graph complex}

In this section we will describe various complexes constructed from the stable ribbon graphs described in the preceding section which form the complex of orbi-cellular chains on moduli spaces of Riemann surfaces. We begin with their definition.

\begin{defi}
The \emph{stable ribbon graph complex}\footnote{Note that in \cite{andrey} the terminology `prestable ribbon graph complex' was used for $\lgc$ and the term `stable ribbon graph complex' was reserved for $\kgc$.} $\lgc$ is the complex whose underlying vector space is freely generated by isomorphism classes of oriented stable ribbon graphs, modulo the relation that reversing the orientation on a stable ribbon graph is equivalent to multiplying by $(-1)$. The differential $\partial$ is given by summing over all possible contractions of the edges:
\[ \partial(\Gamma):=\sum_{e\in E(\Gamma)} \Gamma/e.\]
Note that some edges cannot be contracted, in which case the corresponding term in the sum is defined to be zero. The grading on this complex is given by counting the number of edges. The homology of this complex will be denoted by $\lgh$.
\end{defi}

Note that this complex has a natural subspace (not a subcomplex) which is generated by those stable ribbon graphs for which \emph{the boundary defect of every vertex is zero}. We make the following definition:

\begin{defi}
We define the complex $\kgc$ to be the complex which is generated by those stable ribbon graphs with \emph{everywhere vanishing boundary defect} and whose differential is uniquely defined by the requirement that the natural projection
\begin{equation} \label{eqn_quot1}
\lgc\to\kgc
\end{equation}
is a morphism of complexes. Its homology will be denoted by $\kgh$.
\end{defi}

Furthermore, this complex has a natural subspace (again, not a subcomplex) which is generated by those stable ribbon graphs for which \emph{the boundary defect and the genus defect vanishes at every vertex} and such that \emph{every vertex is partitioned into just one cycle}. In fact, these graphs are just \emph{ribbon graphs}, hence we make the following definition:

\begin{defi}
We define the complex $\gc$ to be the complex which is generated by ribbon graphs and whose differential is uniquely defined by the requirement that the natural projection
\begin{equation} \label{eqn_quot2}
\kgc\to\gc
\end{equation}
be a morphism of complexes. We call this object the \emph{ribbon graph complex}. Its homology will be denoted by $\gh$.
\end{defi}

Note that all three complexes have a natural commutative multiplication given by taking the disjoint union of graphs. Since any graph can be uniquely decomposed into its connected components, it follows that these complexes are freely generated, as algebras, by \emph{connected} graphs; hence they have the canonical structure of differential graded commutative cocommutative Hopf algebras in which the subcomplex generated by connected graphs coincides with the subspace of primitive elements of the Hopf algebra structure. We will denote the homology of this subcomplex by adjoining the prefix $P$.

We can now formulate the fundamental theorem; which follows from the work of Harer, Mumford, Penner, Thurston, Kontsevich and Looijenga; that describes the relationship between these complexes and moduli spaces of Riemann surfaces. For a given locally compact topological space $X$, let us denote the one-point compactification of this space by $X^\infty$. We begin by collecting all the moduli spaces of different genera and with varying numbers of marked points into one object by making the following definitions:
\begin{displaymath}
\begin{split}
\msup & := \bigsqcup_{\begin{subarray}{c} n\geq 1 \\ g>1-\frac{n}{2} \end{subarray}} [\mspc \times \Delta^{\circ}_{n-1}]^{\infty}/\mathbb{S}_n, \\
\ksup & := \bigsqcup_{\begin{subarray}{c} n\geq 1 \\ g>1-\frac{n}{2} \end{subarray}} [\kcmp \times \Delta^{\circ}_{n-1}]^{\infty}/\mathbb{S}_n, \\
\lsup & := \bigsqcup_{\begin{subarray}{c} n\geq 1 \\ g>1-\frac{n}{2} \end{subarray}} \lcmp/\mathbb{S}_n; \\
\end{split}
\end{displaymath}

where the symmetric group $\mathbb{S}_n$ acts naturally on these moduli spaces by the diagonal action which permutes the labels of the marked points and the barycentric coordinates of $\Delta_{n-1}$.

\begin{theorem} \label{thm_decompose}
For all $k\geq 1$ there exists the following commutative diagram
\[ \xymatrix{ H_k\lsup \ar[r] \ar@{=}[d] & H_k\ksup \ar[r] \ar@{=}[d] & H_k\msup \ar@{=}[d] \\ P\lgh[k] \ar[r] & P\kgh[k] \ar[r] & P\gh[k] } \]
\end{theorem}

\begin{proof}
The vertical isomorphisms are just the formal expressions of theorems \ref{thm_orbimod}, \ref{thm_orbikmod} and \ref{thm_orbilmod} that the subcomplexes of $\lgc$, $\kgc$ and $\gc$ generated by connected graphs are precisely the complex of orbi-cellular chains on the corresponding moduli spaces and hence compute the homology of these moduli spaces. The lower horizontal maps are just those induced by \eqref{eqn_quot1} and \eqref{eqn_quot2}.

The top left horizontal map is the morphism which is induced by the mapping
\[ \dmcmp\times\Delta_{n-1} \to  [\kcmp\times\Delta^{\circ}_{n-1}]^\infty \]
which sends any point in $\dmcmp\times\partial\Delta_{n-1}$ to the point at infinity.

Likewise the top right horizontal map is induced by the mapping
\[ [\dmcmp\times\Delta^{\circ}_{n-1}]^{\infty} \to [\mspc\times\Delta^{\circ}_{n-1}]^{\infty} \]
which maps every point in $[\partial\dmcmp\times\Delta^{\circ}_{n-1}]^{\infty}$ to the point at infinity.
\end{proof}

\begin{rem}
Note that standard arguments from algebraic topology allow one to relate the homology of the one-point compactification of the moduli space $[\mspc \times \Delta^{\circ}_{n-1}]/\mathbb{S}_n$ to the homology of the one-point compactification of $\mspc$. Briefly, since the one-point compactification of a cartesian product of two spaces is the smash product of the one-point compactifications of each individual space, the corresponding long exact sequence in homology, which obviously splits, shows that the homology of the one-point compactification of $\mspc \times \Delta^{\circ}_{n-1}$ is just a $n-1$-shifted copy of $H_\bullet\mspc$. Since we are working rationally, taking $\mathbb{S}_n$-coinvariants commutes with homology. Similar arguments apply to $\kcmp \times \Delta^{\circ}_{n-1}$. Note that since $\lcmp$ is not fibered over the moduli space, no such analogue holds for this space.
\end{rem}

\begin{rem}
Given any Riemann surface we can recover the number of marked points $n$ and the arithmetic genus $g$ of this Riemann surface from its (stable) ribbon graph $\Gamma$ using the formulae:
\begin{displaymath}
\begin{split}
n & = n_p + \sum_{v\in V(\Gamma)} n(v), \\
g & = 1 - |V(\Gamma)| + \frac{1}{2}(|E(\Gamma)| + |C(\Gamma)| - n_p) + \sum_{v\in V(\Gamma)} g(v);
\end{split}
\end{displaymath}
where $n_p$ is the number of perimeters of $\Gamma$, $|C(\Gamma)|$ is the total number of cycles of $\Gamma$ and $n(v)$ and $g(v)$ denote the boundary and genus defect respectively.

It follows from this fact that the graph complexes $\lgc$, $\kgc$ and $\gc$ all split as a sum of graph complexes which each have a fixed genera and number of marked points, hence we can equate the homology of each individual summand directly to the homology of the corresponding moduli space having the same genera and number of marked points. However, collecting all these spaces into one object will allow us to provide a convenient formulation of our main theorem in the next section.
\end{rem}

\section{The main theorem} \label{sec_main}

In this section we formulate and prove our main theorem which states that the homology of the moduli space $\lcmp$ is identical to the (stable) relative homology of the differential graded Lie algebra defined by Theorem \ref{thm_dgla}. We then relate the diagram of differential graded Lie algebras defined by \eqref{eqn_seqdglas} to the natural maps on the moduli spaces of Riemann surfaces given by collapsing orbi-cells lying on the boundary to a point. We remind the reader that we \emph{always work over} $\gf$; in particular we view all our differential graded Lie algebras as $\gf$-vector spaces, notwithstanding the fact that they may be naturally defined over larger polynomial rings.

A necessary ingredient in the proof of the main theorem will be the invariant theory for the Lie algebra $\mathfrak{pe}[\mathbb{Q}^{n|n}]$ of linear symplectic vector fields, therefore we begin by recalling the results of \cite{sergeev}.

\begin{defi}
A \emph{chord diagram} is a partition of the set $\{1,\ldots,2k\}$ into pairs. For a given positive integer $k$, we denote the set of all such chord diagrams by $\mathcal{C}(k)$.
\end{defi}

\begin{defi} \label{def_invariant}
Let $(V,\omega)$ be a symplectic vector space. For any chord diagram
\begin{equation} \label{eqn_chord}
c:=\{i_1,j_1\},\ldots,\{i_k,j_k\}
\end{equation}
we can define a linear map
\[ \omega_c:V^{\otimes 2k}\to\gf \]
by the formula
\[ \omega_c(x_1\otimes\cdots\otimes x_{2k}):= (-1)^p\langle x_{i_1},x_{j_1} \rangle\cdots\langle x_{i_k},x_{j_k} \rangle, \]
where $(-1)^p$ is the sign coming from the Koszul sign rule for the permutation
\[ x_1,x_2,\ldots,x_{2k-1},x_{2k} \mapsto x_{i_1},x_{j_1},\ldots,x_{i_k},x_{j_k}.\]
\end{defi}

Note that since the inner product is odd, this sign depends on how the pairs in \eqref{eqn_chord} are ordered. We can get around this issue by assuming that $i_1<i_2<\ldots<i_k$. These maps $\omega_c$ are invariant under the action of the Lie algebra $\mathfrak{pe}[V]$. In fact the following theorem due to Sergeev \cite{sergeev} tells us that they form a basis for all the invariants.

\begin{theorem} \label{thm_invariant}
Let $(V,\omega)$ be a symplectic vector space of dimension $n|n$:
\begin{enumerate}
\item
The set
\[\{\omega_c:V^{\otimes 2k}\to\gf; \ c\in\mathcal{C}(k)\}\]
forms a basis for the space of $\mathfrak{pe}[V]$-invariant linear functions on $V^{\otimes 2k}$, providing that $n\geq k$.
\item
The dimension of the space of $\mathfrak{pe}[V]$-invariant linear functions on $V^{\otimes 2k-1}$ is zero for all $k$.
\end{enumerate}
\end{theorem}
\noproof

For every positive integer $n$ there is a canonical symplectic vector space
\[ \gf^{n|n} :=  \langle \underbrace{x_1,\ldots,x_n}_{even};\underbrace{\xi_1,\ldots,\xi_n}_{odd} \rangle \]
with symplectic form
\[ \omega:=\sum_{i=1}^n dx_i d\xi_i \]
to which any other symplectic vector space of the same dimension is isomorphic.

Let us define new differential graded Lie algebras by taking the stable limit of those defined in Section \ref{sec_dgla}:
\begin{displaymath}
\begin{split}
\Lambda_{\gamma,\nu} & :=\dilim{n}{\left[\Lambda_{\gamma,\nu}[\gf^{n|n}]\right]}, \\
\Lambda_{\gamma} & :=\dilim{n}{\left[\Lambda_{\gamma}[\gf^{n|n}]\right]}, \\
\mathfrak{g} & :=\dilim{n}{\left[\mathfrak{g}[\gf^{n|n}]\right]}, \\
\mathfrak{pe} & :=\dilim{n}{\left[\mathfrak{pe}[\mathbb{Q}^{n|n}]\right]}.
\end{split}
\end{displaymath}

\begin{rem}
The (relative mod $\mathfrak{pe}$) Chevalley-Eilenberg complex of the Lie algebra $\mathfrak{g}$ has a natural commutative multiplication induced by the morphism of Lie algebras
\[ \mathfrak{g}[\mathbb{Q}^{n|n}] \oplus \mathfrak{g}[\mathbb{Q}^{m|m}] \to \mathfrak{g}[\mathbb{Q}^{n+m|n+m}]. \]
Combining this with the usual diagonal on this complex yields the structure of a commutative cocommutative Hopf algebra. Precisely the same remarks apply to the (relative) Chevalley-Eilenberg complexes of the differential graded Lie algebras $\Lambda_{\gamma,\nu}$ and $\Lambda_{\gamma}$.
\end{rem}

Now we introduce a map which formally resembles Wick's formula for integrating with respect to a Gaussian measure.

\begin{defi}
We define a map
\[ I:\rce{\Lambda_{\gamma,\nu}}{\mathfrak{pe}}\to\lgc \]
as follows. A typical element $x$ of $\rce{\Lambda_{\gamma,\nu}}{\mathfrak{pe}}$ is represented by a product
\[ x:=x_1\cdot x_2\cdots x_m \]
of elements $x_i\in\Lambda_{\gamma,\nu}[\gf^{d|d}]$, for some $d>0$. In turn every element $x_i$ is represented by a product
\[ x_i:=\gamma^{g_i}\nu^{n_i}\cdot y^i_1\cdot y^i_2\cdots y^i_{k_i} \]
of elements $y^i_j\in\mathfrak{h_{\geq 1}}[\gf^{d|d}]$ and powers of the deformation parameters $\gamma$ and $\nu$. Finally, each element $y^i_j$ is represented by a product
\[ y^i_j:=z^i_{j1}\cdot z^i_{j2}\cdots z^i_{jl_{ij}} \]
of elements $z^i_{jr}\in(\gf^{d|d})^*$. Hence the total number of tensors we have is the sum of the $l_{ij}$. If the total number of tensors are odd then we define $I(x)$ to be zero, hence we assume that $\sum_{i,j} l_{ij}=2M$ is even.

Now, for each chord diagram
\[ c:=\{i_1,j_1\},\ldots,\{i_M,j_M\} \]
there is an obvious way to construct a corresponding graph. Namely, we take a graph having $m$ vertices and partition the $i$th vertex into $k_i$ cycles, such that the $j$th cycle has valency $l_{ij}$. We set the genus defect at this vertex to be $g_i$ and the boundary defect to be $n_i$. The chord diagram $c$ provides a way to pair up the half-edges of the graph. The orientation can be defined canonically by assuming that $i_1<\ldots <i_M$. Let us denote this graph by $\Gamma_c$.

Finally, the map $I$ is defined by the formula:
\[ I(x):=\sum_{c\in\mathcal{C}(M)} \omega_c(x)\Gamma_c. \]
That is to say that the coefficient of a graph $\Gamma$ is determined by first placing the tensors $x_i$ at the vertices of $\Gamma$, with the subtensors $y^i_j$ placed on the cycles of that vertex using the cyclic ordering, then contracting these tensors by applying the inner product $\innprod^{-1}$ to each edge.
\end{defi}

The map $I$ can obviously be restricted to the subspaces (not subcomplexes) $\rce{\Lambda_\gamma}{\mathfrak{pe}}$ and $\rce{\mathfrak{g}}{\mathfrak{pe}}$. This leads to the following commutative diagram:
\begin{equation} \label{eqn_maindiagram}
\xymatrix{\lgc \ar[r] & \kgc \ar[r] & \gc \\ \rce{\Lambda_{\gamma,\nu}}{\mathfrak{pe}} \ar[r] \ar[u]^I & \rce{\Lambda_\gamma}{\mathfrak{pe}} \ar[r] \ar[u]^I & \rce{\mathfrak{g}}{\mathfrak{pe}}\ar[u]^I}
\end{equation}
where the top horizontal maps are those defined by \eqref{eqn_quot1} and \eqref{eqn_quot2} and the bottom horizontal maps are those defined by diagram \eqref{eqn_seqdglas}.

Now we are ready to formulate the main theorem.

\begin{theorem} \label{thm_main}
The vertical maps of diagram \eqref{eqn_maindiagram} are isomorphisms of differential graded Hopf algebras.
\end{theorem}

We have as an immediate corollary:

\begin{cor}
Combining Theorem \ref{thm_decompose} with Theorem \ref{thm_main} yields the following commutative diagram for all $k\geq 1$ relating the homology of the differential graded Lie algebras defined in Section \ref{sec_constructdgla} to the homology of the compactifications of the moduli space defined in Section \ref{sec_js}:
\[ \xymatrix{H_k\lsup \ar[r] & H_k\ksup \ar[r] & H_k\msup \\ P\hrce[k]{\Lambda_{\gamma,\nu}}{\mathfrak{pe}} \ar[r] \ar[u]^I & P\hrce[k]{\Lambda_\gamma}{\mathfrak{pe}} \ar[r] \ar[u]^I & P\hrce[k]{\mathfrak{g}}{\mathfrak{pe}}\ar[u]^I} \]
\end{cor}

\begin{rem}
It is in fact possible to consider other compactifications of the moduli space of curves fitting into the above diagram which correspond to setting one or both of the deformation parameters $\gamma$ and $\nu$ to zero. However, in order to avoid unduly complicating the exposition of this paper, this perspective will not be pursued.
\end{rem}

\begin{rem}
It is very likely that Theorem \ref{thm_main} could be generalised to the setting of an arbitrary modular operad. Given any modular operad, one can associate to it a certain graph complex by decorating the vertices of the graphs by this modular operad. One should also be able to associate a differential graded Lie algebra to this modular operad which recovers the homology of this graph complex. For a treatment of graph complexes from the perspective of modular operads, the reader may consult \cite{andrey}.
\end{rem}

\begin{proof}[Proof of Theorem \ref{thm_main}]
First we explain why the vertical maps $I$ must be bijective. This is just a direct application of Theorem \ref{thm_invariant} which describes the invariants of the Lie algebra $\mathfrak{pe[\mathbb{Q}^{n|n}]}$. We can construct an explicit inverse to $I$ as follows. To each graph $\Gamma$ we define a certain tensor $x_\Gamma$ in the Chevalley-Eilenberg complex $\rce{\Lambda_{\gamma,\nu}[\mathbb{Q}^{n|n}]}{\mathfrak{pe}[\mathbb{Q}^{n|n}]}$, where $n$ is the number of edges of the graph. We do this by decorating every edge of the graph with a pair of tensors $x_i,\xi_i$ for $i=1,\ldots, n$ so that every edge of the graph is decorated with a distinct pair of tensors. The structure of the vertices of $\Gamma$ and their cycles gives us an obvious way to interpret this object as a tensor $x_\Gamma$ in the Chevalley-Eilenberg complex $\rce{\Lambda_{\gamma,\nu}[\mathbb{Q}^{n|n}]}{\mathfrak{pe}[\mathbb{Q}^{n|n}]}$. Since the tensors defined in this manor are just $\mathfrak{pe}[\mathbb{Q}^{n|n}]$-coinvariants which are obviously dual to the $\mathfrak{pe}[\mathbb{Q}^{n|n}]$-invariants $\omega_c$ described in Definition \ref{def_invariant}, applying the map $I$ to the tensor $x_\Gamma$ will give us back the graph $\Gamma$. By applying Theorem \ref{thm_invariant}, we may assume that any tensor in $\rce{\Lambda_{\gamma,\nu}}{\mathfrak{pe}}$ can be represented by one of the form $x_\Gamma$ for some graph $\Gamma$, hence the map $I$ is bijective.

Next, we mention that an essentially standard calculation verifies that the map $I$ is a map of Hopf algebras. Since the proof of this fact involves essentially the same argument as that employed in Theorem 4.18 of \cite{char}, we choose not repeat it here.

Thus, it remains to check that the map $I$ is a map of complexes. This is technically much simpler if we assume that all our tensors have the form $x_\Gamma$ for some graph $\Gamma$, which we may do as a consequence of the results on the invariant theory for $\mathfrak{pe}[\mathbb{Q}^{n|n}]$ as explained above. A straightforward calculation using formula \eqref{eqn_bvbracket} for the Lie bracket on $\Lambda_{\gamma,\nu}$, formula \eqref{eqn_bracket} for the Lie bracket on noncommutative 0-forms and formula \eqref{eqn_cobracket} for the Lie cobracket on $0$-forms then verifies that $I$ is a map of complexes. Intuitively, the correspondence is clear:

\begin{enumerate}
\item
The term in the differential on $\rce{\Lambda_{\gamma,\nu}}{\mathfrak{pe}}$ contributed by the Lie bracket on $\Lambda_{\gamma,\nu}$ corresponds to contracting edges of the form described in Definition \ref{def_contract} \eqref{item_contractedge}:
\begin{figure}[htp]
\centering
\includegraphics{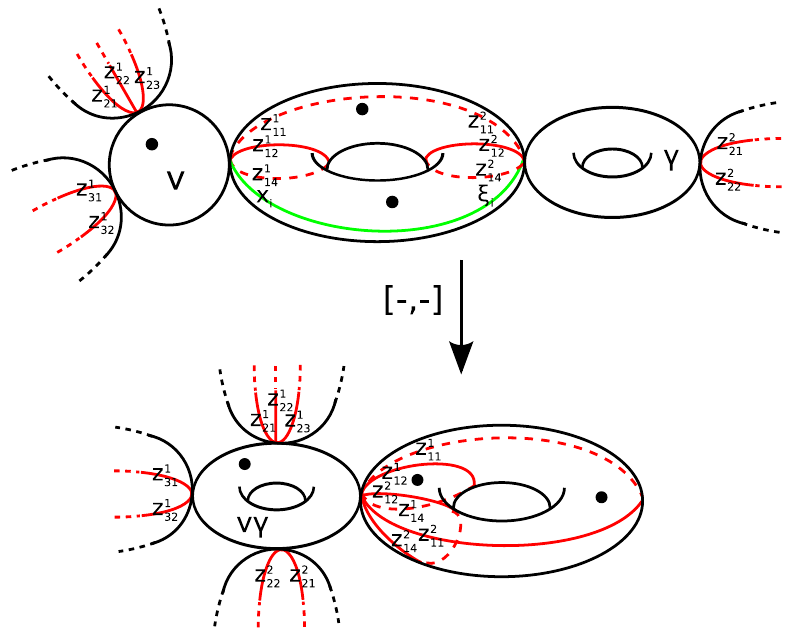}
\caption{The contribution from the bracket on $\Lambda_{\gamma,\nu}$.}
\end{figure}
\item
The term in the differential contributed by the Lie bracket on 0-forms corresponds to contracting edges of the form described in Definition \ref{def_contract} \eqref{item_contractloop1}:
\begin{figure}[htp]
\centering
\includegraphics{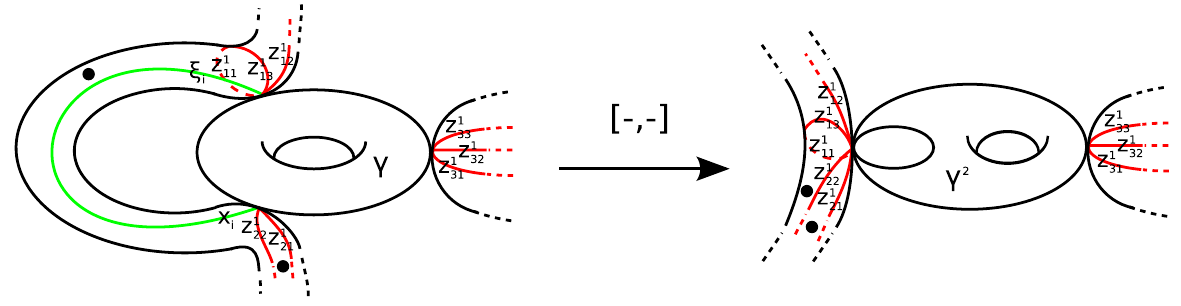}
\caption{The contribution from the bracket on $\mathfrak{h}_{\geq 1}$.}
\end{figure}
\item
The term in the differential contributed by the Lie cobracket on 0-forms corresponds to contracting edges of the form described in Definition \ref{def_contract} \eqref{item_contractloop2}:
\begin{figure}[htp]
\centering
\includegraphics{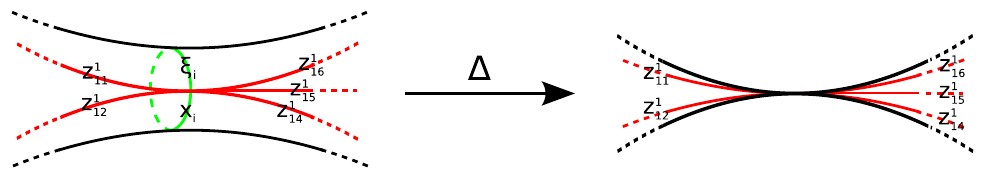}
\caption{The contribution from the cobracket on $\mathfrak{h}_{\geq 1}$.}
\end{figure}
\end{enumerate}
\end{proof}

\end{document}